\documentclass[12pt]{article}
\usepackage{amssymb}
\usepackage{color}
\usepackage{tikz}
\newcommand{\bK}{{\mathbf K}}
\newtheorem{teorema}{Theorem}[section]
\newtheorem{lema}[teorema]{Lemma}
\newtheorem{nota}[teorema]{Remark}
\newtheorem{ejemplo}[teorema]{Example}
\newcommand{\NWP}{\mathrm{gcd}}
\newcommand{\ord}{\mathrm{ord}}

\newenvironment{proof}[1][Proof]{\textbf{#1.} }{\
\rule{0.5em}{0.5em}}

\begin{document}
\title{ Euclidean algorithm and polynomial equations after Labatie
\footnotetext{
     \noindent   \begin{minipage}[t]{5in}
       {\small
       2000 {\it Mathematics Subject Classification:\/}Primary 12;
       Secondary 14H20.\\
       Key words and phrases: polynomial equations,
       Euclidean algorithm, intersection multiplicity.\\
       The first-named author was partially supported by the Spanish Project
       PNMTM 2007-64007.}
       \end{minipage}}}
\author{E.R. Garc\'{\i}a Barroso and A. P\l{}oski}
\maketitle
\begin{abstract}
\noindent We recall Labatie's effective method
of solving polynomial equations with two unknowns
by using the Euclidean algorithm.
\end{abstract}

\section*{Introduction}

\noindent The French mathematician Labatie published in 1835 a booklet on a
method of solving polynomial systems of equations in two unknowns (see \cite{Finck1}). He used the polynomial division
 to replace the given system of equations by the collection of triangular systems. Labatie's
theorem can be found in some old Algebra books: by Finck \cite{Finck2}, Serret \cite{Serret} and Netto \cite{Netto}, but as far as we know, not in any Algebra text book written in the twentieth century.

\medskip

\noindent In this paper we recall Labatie's method following Serret \cite{Serret} (pp. 196-206). Then we give, in a modern setting, an improvement of Labatie's result due to Bonnet \cite{Bonnet}.

\medskip

\noindent Let $\mathbf K$ be a field of arbitrary characteristic. We shall consider polynomials with coefficients in $\bK$.
If  $W=W(x,y)\in\bK[x,y]$ then we denote by $\deg_yW$ the degree of $W$ with respect to
$y$. We say that a non-zero polynomial $W$ is $y${\em -primitive} if it is a primitive polynomial in the ring
$\bK[x][y]$, that is, if $1$ is the greatest common divisor of all the non-zero coefficients that are dependent on $x$.
If  $V$, $W\in\bK[x,y]$ satisfy the condition $0<\deg_yV\leq\deg_yW$ then there are polynomials
$Q$ (quotient), $R$ (remainder) in $\bK[x,y]$ and a non-zero polynomial $u=u(x)\in\bK[x]$ such that  $uW=QV+R$, where $\deg_yR<\deg_yV$
or $R=0$

\medskip

\noindent The greatest common divisor of polynomials $V$, $W$ may be computed  using the Euclidean algorithm, see
\cite{Bocher} chapter XVI. Recently Hilmar and Smyth~\cite{Hilmar-Smyth} gave a very simple proof of B\'ezout's theorem
for plane projective curves using as a main tool the Euclidean division.


\section{Euclidean algorithm}
Let $V_1$, $V_2\in\bK[x,y]$ be coprime and $y$-primitive polynomials such that $0<\deg_yV_2\leq\deg_yV_1$.

\medskip

\noindent Using the polynomial division we get a sequence of $y$-primitive polynomials
$V_3$, \dots, $V_{n+1}$ of decreasing $y$-degrees $0<\deg_yV_{n+1}<\dots<\deg_yV_3< deg_yV_2$
such that
\begin{eqnarray*}
   u_1V_1 &=& Q_1V_2+v_1V_3, \\
   u_2V_2 &=& Q_2V_3+v_2V_4, \\
   &\vdots& \\
   u_{n-1}V_{n-1} &=& Q_{n-1}V_n+v_{n-1}V_{n+1}, \\
   u_nV_n &=& Q_nV_{n+1}+v_n, \\
\end{eqnarray*}
where $u_1$,\dots, $u_n$, $v_1$,\dots, $v_n$ are non-zero polynomials of the ring $\bK[x]$.
Let be $V_{n+2}=1$ and write the above equalities in the form
$$
u_{i}V_{i}=Q_{i}V_{i+1}+v_{i}V_{i+2}\quad\mbox{for $i=1,\dots,n$.}
\leqno (1)_{i}
$$

\noindent In what follows we call $n$ the number of steps performed by the Euclidean algorithm on
input $(V_1$, $V_2)$. We will keep the above notation in all this note.

\section{Labatie's elimination}

Let us define two sequences $d_1,\ldots,d_n$ and $w_1,\ldots,w_n$ of polynomials in $x$ determined by the sequences
$u_1,\ldots,u_n$ and $v_1,\ldots,v_n$ in a recurrent way. We let
$d_1=\NWP(u_1,v_1)$, $w_1=\frac{u_1}{d_1}$ and
 $d_{i}=\NWP(w_{i-1}u_{i},v_{i})$,  $w_{i}=\frac{w_{i-1}u_{i}}{d_{i}}$
for $i\in \{2,\ldots,n\}$. It is easy to see that $w_{i}=
\frac{u_1\cdots u_{i}}{d_1\cdots d_{i}}$ in $\bK[x]$  for all $i \in \{1,\ldots,n\}$.

\medskip

\noindent For any $V$, $W\in\bK[x,y]$ we let
$\{V=0,W=0\}=\{\,P\in\bK^2: V(P)=W(P)=0\,\}$.

\begin{teorema}[Labatie 1835]\label{Labatie}
With notations and assumptions given above we have
$$ \left\{V_1=0,V_2=0\right\}=\bigcup_{i=1}^n \left\{V_{i+1}=0,\frac{v_{i}}{d_{i}}=0\right\} .
$$
\end{teorema}

\noindent We present the proof of the above theorem in Section \ref{proof-Labatie}.

\noindent Labatie's theorem shows that the system of equations $V_1(x,y)=0$, $V_2(x,y)=0$ is equivalent
to the collection of triangular systems

$$ V_{i+1}(x,y)=0,\; \frac{v_{i}}{d_{i}}(x)=0\qquad (i=1,\dots,n) .
$$

\noindent  Labatie's theorem fell into oblivion for a long time. At the beginning of the 1990's Lazard in 
\cite{Lazard} proved that every system of polynomial equations in many unknowns with a finite
number of solutions in the algebraic closure of $\bK$ is equivalent to the union of triangular systems, which can be obtained from Gr\"obner bases. Kalkbrener in \cite{Kalkbrener1} and  \cite{Kalkbrener2}
developed the theory of elimination sequences based on the Euclidean algorithm. His method of computing solutions of systems of polynomials equations turned out to be very efficient if applied to 
systems of two or three unknowns (see  \cite{Kalkbrener2} and the references given therein for the comparison with Gr\"obner basis methods). Neither Lazard nor Kalkbrener mentioned Labatie's work. Only Glashof in \cite{Glashoff} recalled Labatie's method after Netto \cite{Netto} and compared it with
Kalkbrener's approach to polynomials equations. In what follows we need the notion of multiplicity of a solution of a system of two equations in two unknowns. The definition we are going to present is quite sophisticated. The reader not acquainted with it may assume the five properties of multiplicity given below as axiomatic definition of this notion.

\medskip

\noindent Let $P\in \bK^2$. We define the local ring of rational functions regular at $P$ to be
$$ \bK[x,y]_P=\left \{\frac{R}{S}\;:\; R,S \in \bK[x,y], S(P)\neq 0\right \}.$$

\noindent The ring $ \bK[x,y]_P$ is a unique factorization domain. The units of $ \bK[x,y]_P$ are rational functions
$\frac{R}{S}$ such that $R(P)S(P)\neq 0$.

\medskip

\noindent Let $(V,W)_P$ be the ideal generated by polynomials $V$ and $W$ in $ \bK[x,y]_P$. Following
\cite{Fulton}, we define the {\em intersection multiplicity} $i_P(V,W)$ to be the dimension of the $\bK$-vector space
$ \bK[x,y]_P/(V,W)_P$. We call also $i_P(V,W)$ the {\em multiplicity of the solution} $P$ of the system $V=0$, $W=0$.

\medskip

\noindent Let us recall the basic properties of the intersection multiplicity which hold for any field $\bK$ (not
necessarily algebraically closed):

\begin{enumerate}
\item $i_P(V,W)<+\infty$ if and only if $P\not\in \{\gcd(V,W)=0\}$,
\item $i_P(V,W)>0$ if and only if $P\in \{V=W=0\}$,
\item $i_P(V,WW')=i_P(V,W)+i_P(V,W')$,
\item $i_P(V,W)$ depends only on the ideal $(V,W)_P$. \\ Intuitively: $i_P(V,W)$ does not change when we replace the system $V=0$, $W=0$ by another one equivalent to it near $P$.

\noindent Moreover, it is easy to check that

\item if $P=(a,b)$ is a solution of the triangular system $W(x,y)=0$, $w(x)=0$ then
$i_P(W,w)=(\ord_a w)(\ord_b W(a,y))$, where $\ord_c p$ denotes the multiplicity of the root
$c$ in the polynomial $p=p(x)\in \bK[x]$. By convention $\ord_c p=0$ if $p(c)\neq 0$.

\noindent The following example may be helpful to acquire an intuition of intersection multiplicity. Let us consider the parabola $y^2-x=0$ over the field of real numbers. Applying Property 5 to the triangular system $y ^2-x=0$, $x-c=0$ we check that the axis $x=0$ intersects the parabola in $(0,0)$ with multiplicity 2 but the line $x-c=0$,
where $c >0$ intersects it in two points $(c,\sqrt{c} )$  and $(c,-\sqrt{c})$, each with multiplicity 1. If
$c\to   0^+$ then the two points coincide.

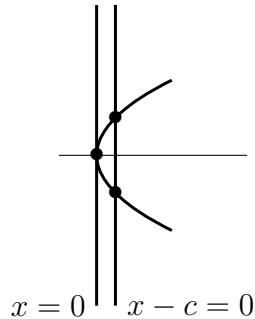
\begin{center}
\begin{tikzpicture}
\begin{scope}[rotate=270]
\draw(-2,0) -- (2,0) node[left] {$x=0$}; %
      \draw(0,-0.5) -- (0,2) ;%
\draw  [domain=-2:2,black] [line width=0.4mm] plot (\x,1/4)node[right] {$x-c=0$};
\draw  [domain=-2:2,black]  [line width=0.4mm] plot (\x,0);
\draw [line width=0.4mm] (-1,1) parabola[parabola height=-1cm] (1,1);
\draw (0,0) node {$\bullet$} ;
\draw (-1/2,1/4) node {$\bullet$} ;
\draw (1/2,1/4) node {$\bullet$} ;
\end{scope}
\end{tikzpicture}
\end{center}

\end{enumerate}

\medskip

\noindent Note also that the system of equations $y^2-x=0$, $x-c=0$ has for $c\neq 0$ two complex solutions, which are arbitrary close to the origin for small enough complex $c$. This observation leads to the  {\em dynamic definition} of intersection multiplicity for algebraic complex curves (see \cite{Teissier}, Section 6).

\medskip

\noindent The following theorem due to Bonnet \cite{Bonnet} is an improvement of Labatie's result:

\begin{teorema}[Bonnet 1847]
\label{Bonnet}
For any $P\in \bK^2$ we have
$$ i_P(V_1,V_2)=\sum_{i=1}^n i_P\Bigl(V_{i+1},\frac{v_{i}}{d_{i}}\Bigr).
$$
\end{teorema}

\noindent Bonnet, like Labatie, considered polynomials with complex coefficients and used the definition
of the intersection multiplicity in terms of Puiseux series. In Section \ref{proof-Bonnet}  we present
a short proof of Theorem \ref{Bonnet} based on Labatie's calculations (Section \ref{auxiliary}) and the properties
of the intersection multiplicity listed above.

\begin{ejemplo}
Let $V_1=y^5-x^3$, $V_2=y^3-x^4$. Using the Euclidean algorithm we get $y^5-x^3=y^2(y^3-x^4)+x^3(xy^2-1)$,
$x(y^3-x^4)=y(xy^2-1)+y-x^5$ and $xy^2-1=(xy+x^6)(y-x^5)+x^{11}-1$. Hence we have $(u_1,u_2,u_3)=(1,x,1)$, $(v_1,v_2,v_3)=(x^3, 1,x^{11}-1)$ and $(d_1,d_2,d_3)=(1,1,1)$. By Labatie's theorem, we get
\begin{eqnarray*}
\{y^5-x^3=0,y^3-x^4=0\}=&&\\
\{y^3-x^4=0,x^3=0\}\cup \{xy^2-1=0,1=0\}\cup \{y-x^5=0,x^{11}-1=0\}.&&
\end{eqnarray*}

\noindent Therefore the systems $V_1=0$, $V_2=0$ has two solutions $(0,0)$ and $(1,1)$ in $\bK$ and ten solutions in the algebraic closure of $\bK$. To compute the multiplicities of the solutions we use Bonnet's theorem:
$$i_0(y^5-x^3,y^3-x^4)=i_0(y^3-x^4,x^3)+i_0(xy^2-1,1)+i_0(y-x^5,x^{11}-1)=3 \cdot 3+0+0=9.$$ 

\noindent The remaining multiplicities are equal to one. Thus the system $V_1=0$, $V_2=0$ has $9+11=20$ solutions counted with multiplicities.
\end{ejemplo}
\section{Auxiliary lemmas}
\label{auxiliary}
Recall that the polynomials
 $w_{i}$ and $\frac{v_{i}}{d_{i}}$ are coprime.

\begin{lema}\label{Lem}
There exist two sequences of polymomials
 $G_0$,\dots, $G_n$ and
$H_0$,\dots, $H_n$ in the ring $\bK[x,y]$  such that
$$
w_{i-1}V_1=G_{i-1}V_{i}+G_{i-2}V_{i+1}\frac{v_{i-1}}{d_{i-1}}, \leqno (2)_{i}
$$
$$
w_{i-1}V_2=H_{i-1}V_{i}+H_{i-2}V_{i+1}\frac{v_{i-1}}{d_{i-1}} \leqno (3)_{i}
$$
for $i\in \{2,\dots, n+1\}$.
\end{lema}

\noindent \begin{proof} We proceed by induction on $i$. Let's check the first identity.
From the equality  $u_1V_1=Q_1V_2+v_1V_3$ it follows that $d_1=\NWP(u_1,v_1)$ divides the product $Q_1V_2$
and consequently the polynomial $Q_1$ since $V_2$ is $y$-primitive.
Letting $G_0=1$, $G_1=\frac{Q_1}{d_1}$ we get $w_1V_1=G_1V_2+G_0V_3\frac{v_1}{d_1}$
that is $(2)_2$.
Suppose now that $2\leq i<n+1$  and that for some polynomials $G_{i-1}$
and~$G_{i-2}$ the identity ~$(2)_{i}$ holds.
Multiplying the identity $(2)_{i}$ by the polynomial  $u_{i}$ we get
$$
 w_{i-1}u_{i}V_1=
 u_{i}G_{i-1}V_{i}+u_{i}G_{i-2}V_{i+1}\frac{v_{i-1}}{d_{i-1}} .
$$

\noindent Let us insert to the identity above $u_{i}V_{i}=Q_{i}V_{i+1}+v_{i}V_{i+2}$.
After simple computations we get:
$$
 w_{i-1}u_{i}V_1=
 \Bigl(G_{i-1}Q_{i}+u_{i}G_{i-2}\frac{v_{i-1}}{d_{i-1}}\Bigr)V_{i+1}+G_{i-1}v_{i}V_{i+2}.
$$
Since $d_{i}=\NWP(w_{i-1}u_{i},v_{i})$ and the polynomial $V_{i+1}$ is $y$-primitive we get that
$G_{i}:=\frac{G_{i-1}Q_{i}}{d_{i}}+G_{i-2}\frac{u_{i}v_{i-1}}{d_{i}d_{i-1}}$
is a polynomial and we have
$$ w_{i}V_{1}=G_{i}V_{i+1}+G_{i-1}V_{i+2}\frac{v_{i}}{d_{i}}, $$
which is the identity $(2)_{i+1}$. This proves the first part of the lemma.

\noindent To prove the identity $(3)_{i}$ note that
$$ w_1V_2=H_1V_2+H_0V_3\frac{v_1}{d_1} $$
if we let $H_0=0$ and $H_1=\frac{u_1}{d_1}$. This proves $(3)_2$.
To check $(3)_{i}$ we proceed analogously to the proof of $(2)_{i}$ : it suffices
to replace $G_{i}$ by $H_{i}$.
\end{proof}

\begin{nota}
The polynomials $G_{i}$ are defined by $G_0=1$, $G_1=\frac{Q_1}{d_1}$,
$G_{i}=\frac{G_{i-1}Q_{i}}{d_{i}}+\frac{G_{i-2}u_{i}v_{i-1}}{d_{i-1}d_{i}}$
and the polynomials $H_{i}$ by $H_0=0$, $H_1=\frac{u_1}{d_1}$ and
$H_{i}=\frac{H_{i-1}Q_{i}}{d_{i}}+\frac{H_{i-2}u_{i}v_{i-1}}{d_{i-1}d_{i}}$.
\end{nota}

\begin{lema}\label{L2}
With the notations of Lemma \ref{Lem} we have the identities
$$
 (-1)^{i}\frac{v_1\cdots v_{i-1}}{d_1\cdots d_{i-1}}
 V_{i+1}=H_{i-1}V_1-G_{i-1}V_2 \qquad \mbox{for $i\in\{2,\dots,n+1\}$} .
\leqno (4)_{i}
$$
\end{lema}

\noindent \begin{proof}
Let $D_{i}=G_{i}H_{i-1}-G_{i-1}H_{i}$ for $i\in\{2,\dots, n\}$.
Consider the system of equations $(2)_{i}$, $(3)_{i}$ as a linear system with unknowns
$V_{i}$, $V_{i+1}\frac{v_{i-1}}{d_{i-1}}$ with determinant equal to  $D_{i-1}$. Using Cramer's rule we get
\begin{eqnarray*}
 D_{i-1}V_{i} &=& w_{i-1}\left(H_{i-2}V_1-G_{i-2}V_2\right) , \\
 D_{i-1}V_{i+1}\frac{v_{i-1}}{d_{i-1}} &=& -w_{i-1}(H_{i-1}V_1-G_{i-1}V_2) . \\
\end{eqnarray*}

\noindent Replacing in the first equality $i$ by $i+1$ we obtain

\begin{equation}
\label{aaa}
D_{i}V_{i+1}=w_{i}(H_{i-1}V_1-G_{i-1}V_2).
\end{equation}

\noindent Multiplying the second equality by $\frac{u_{i}}{d_{i}}$ we get

\begin{equation}
\label{bbb}
D_{i-1}V_{i+1}\frac{v_{i-1}}{d_{i-1}}\frac{u_{i}}{d_{i}}=
 -w_{i}(H_{i-1}V_1-G_{i-1}V_2).
\end{equation}

\noindent Comparing the left sides of $(\ref{aaa})$ and $(\ref{bbb})$ and cancelling
$V_{i+1}$ we have $D_{i}=-\frac{v_{i-1}u_{i}}{d_{i-1}d_{i}}D_{i-1}$.
Moreover $D_1=G_1H_0-G_0H_1=-\frac{u_1}{d_1}$ and by  induction we have
$$ D_{i}=(-1)^{i}w_{i}\frac{v_1\cdots v_{i-1}}{d_1\cdots d_{i-1}}
$$
which inserted into formula $(\ref{aaa})$
gives the identity $(4)_{i}$.
\end{proof}

\section{Proof of  Labatie's theorem}
\label{proof-Labatie}
\medskip\noindent
We can now give  the proof of Theorem \ref{Labatie}: fix a point
 $P\in\bK^2$. If  $V_{i}(P)=\frac{v_{i-1}}{d_{i-1}}(P)=0$
for a value $i \in\{2,\dots,n+1\}$ then from Lemma
\ref{Lem} it follows that
$V_1(P)=V_2(P)=0$ given that $w_{i-1}(P)\neq0$
since $w_{i-1}$, $\frac{v_{i-1}}{d_{i-1}}$ are coprime.

\noindent Suppose now that $V_1(P)=V_2(P)=0$.
From the identity $(4)_{n+1}$ of Lemma \ref{L2} we get  $\frac{v_1\cdots v_n}{d_1\cdots d_n}(P)=0$. Therefore at least one
of polynomials $\frac{v_1}{d_1}$,\dots, $\frac{v_n}{d_n}$ vanishes at $P$.
If $\frac{v_1}{d_1}(P)=0$ then $P\in\{V_2=\frac{v_1}{d_1}=0\}$.

\noindent If the smallest index  $i$ for which $\frac{v_{i}}{d_{i}}(P)=0$
is strictly greater than 1 then we get, by the identity ~$(4)_{i}$ ,
that $V_{i+1}(P)=0$ because $\frac{v_1\cdots v_{i-1}}{d_1\cdots d_{i-1}}(P)\neq0$ by the definition of $i$.
This proves the theorem.

\section{Proof of  Bonnet's theorem}
\label{proof-Bonnet}
Fix a point $P\in \bK^2$. If $\frac{v_1\cdots v_{n}}{d_1\cdots d_{n}}(P)\neq0$ then by $(4)_{n+1}$ we get
\begin{equation}
\label{ccc} 1\in (V_1,V_2)_P
\end{equation}

\noindent which implies $i_P(V_1,V_2)=0$.

\noindent On the other hand we have $i_P\left(V_{i+1},\frac{v_{i}}{d_{i}}\right)=0$ since
$\frac{v_{i}}{d_{i}}(P)\neq 0$ for $i \in\{1,\ldots,n\}$ and the theorem holds in the case under consideration.

\medskip

\noindent Suppose now that $\frac{v_1\cdots v_{n}}{d_1\cdots d_{n}}(P)=0$ and let $i_0$ be the
smallest index $i \in \{1,\ldots,n\}$ such that $\frac{v_{i_0}}{d_{i_0}}(P)=0$. Therefore we have
$w_{i_0}(P)\neq 0$ since  $\frac{v_{i_0}}{d_{i_0}}$ and $w_{i_0}$ are coprime. Let us check
that

\begin{equation}
\label{ddd}
(V_1,V_2)_P=\left(V_{i_0+1},V_{i_0+2}\frac{v_{i_0}}{d_{i_0}}\right)_P. 
\end{equation}

\noindent From $(2)_{i_0+1}$ and  $(3)_{i_0+1}$ we get

\begin{equation}
\label{zzz}
V_1,V_2\in \left(V_{i_0+1},V_{i_0+2}\frac{v_{i_0}}{d_{i_0}}\right)_P. 
\end{equation}

\noindent On the other hand,  from $(4)_{i_0}$ (if $i_0>1$, the case $i_0=1$ being obvious), we obtain
\begin{equation}
\label{xxx}
V_{i_0+1}\in (V_1,V_2)_P 
\end{equation}

\noindent and  from $(4)_{i_0+1}$, we have
\begin{equation}
\label{yyy}
\frac{v_{i_0}}{d_{i_0}}V_{i_0+2}\in (V_1,V_2)_P. 
\end{equation}

\noindent Combining $(\ref{zzz})$, $(\ref{xxx})$ and $(\ref{yyy})$ we get $(\ref{ddd})$. Equality $(\ref{ddd})$ and the additive property of intersection multiplicity
give

\begin{equation}
\label{ppp}
i_P(V_1,V_2)=i_P\left(V_{i_0+1},\frac{v_{i_0}}{d_{i_0}}\right)+i_P(V_{i_0+1},V_{i_0+2}). 
\end{equation}

\noindent If $i_0=n$ then $(\ref{ppp})$  reduces to 

\begin{equation}
\label{fff}
i_P(V_1,V_2)=i_P\left(V_{n+1},\frac{v_n}{d_n}\right)
\end{equation}

\noindent since $V_{n+2}=1$.

\medskip

\noindent To prove Theorem \ref {Bonnet} we shall proceed by induction on the number $n$ of steps performed by the Euclidean
algorithm. For $n=1$ the theorem follows from $(\ref{fff})$ since $n=1$ implies $i_0=1$. Let $n>1$ and suppose that the theorem holds for all pairs of polynomials for which the number of steps performed by
the Euclidean algorithm is strictly less than $n$.\\

\medskip

\noindent We assume that $i_0<n$ since for $i_0=n$ the theorem is true by $(\ref{fff})$.

\medskip

\noindent Let us put $\overline{V}_{j}=V_{i_0+j}$, where
$j\in \{1,2,\ldots,n-i_0+2\}$. The number of steps performed by the Euclidean algorithm on input $(\overline{V}_1$,
$\overline{V}_2)$ is equal to $\overline{n}=n-i_0<n$. We have $\overline{u}_{j}=u_{i_0+j}$
and $\overline{v}_{j}=v_{i_0+j}$ for $j\in\{1,\ldots,\overline{n}\}$. To relate 
 $\overline{d}_{j}$ and $d_{i_0+j}$ we introduce some notation. 
We will write $u\sim \tilde{u}$ for 
 polynomials $u,\tilde{u}$ associated in the local ring $\bK[x,y]_P$. If $u,\tilde{u} \in \bK[x]$ then
$u\sim \tilde{u}$ if and only if there exist polynomials $r,s\in  \bK[x]$ such that $su=r\tilde{u}$ and $r(P)s(P)\neq 0$. Note that $\gcd(u,v)\sim \gcd(\tilde{u},v) $ if $u\sim \tilde{u}$. We claim that

\begin{equation}
\label{eee}
\overline{d}_{j}\sim d_{i_0+j},\;\;\overline{w}_{j}\sim w_{i_0+j}\;\;\hbox{\rm for }\; j\in\{1,\ldots,\overline{n}\}.
\end{equation}

\noindent Let us check (\ref{eee}) by induction on $j$.

\medskip

\noindent If  $j=1$ then $\overline{d}_{{1}}=\gcd(\overline{u}_1,\overline{v}_1)=\gcd(u_{i_0+1},v_{i_0+1})
\sim \gcd(w_{i_0}u_{i_0+1},v_{i_0+1})=d_{i_0+1}$ since $w_{i_0}\sim 1$. Hence we get $\overline{w}_1=\frac{\overline{u}_1}{\overline{d}_1}=\frac{u_{i_0+1}}{\overline{d}_1}\sim \frac{w_{i_0}u_{i_0+1}}{d_{i_0+1}}$, which proves $(\ref{eee})$ for $j=1$.

\medskip

\noindent Suppose that $(\ref{eee})$ holds for a $j<\overline{n}$. Then we get $$\overline{d}_{j+1}=\gcd(\overline{w}_{j}\overline{u}_{j+1},\overline{v}_{j+1})\sim \gcd(w_{i_0+j}u_{i_0+j+1},v_{i_0+j+1})=d_{i_0+j+1}$$

\noindent since $\overline{w}_{j}\sim w_{i_0+j}$ by the inductive assumption, and
$$\overline{w}_{j+1}=\frac{\overline{w}_{j}\overline{u}_{j+1}}{\overline{d}_{j+1}}\sim \frac{w_{i_0+j}u_{i_0+j+1}}{d_{i_0+j+1}}=w_{i_0+j+1}.$$

\noindent This finishes the proof of $(\ref{eee})$.

\medskip

\noindent Now we can finish the proof of the theorem. By the inductive assumption applied to the pair $\overline{V}_1$,
$\overline{V}_2$ we get

\begin{eqnarray*}
i_P(V_{i_0+1},V_{i_0+2})&=&i_P(\overline{V}_1,
\overline{V}_2)=\sum_{j=1}^{\overline{n}}i_P\left(\overline{V}_{j+1},
\frac{\overline{v}_{j}}{\overline{d}_{j}}\right)\\
&=&
\sum_{j=1}^{\overline{n}}i_P\left(V_{i_0+j+1},
\frac{v_{i_0+j}}{d_{i_0+j}}\right)=\sum_{i=i_0+1}^n i_P\left(V_{i+1},\frac{v_{i}}{d_{i}}\right)\\
\end{eqnarray*}

\noindent since $\overline{d}_{j}\sim d_{i_0+j}$ by $(\ref{eee})$ which together with $(\ref{ppp})$
 proves the inductive step and so the theorem.



\medskip
\noindent
{\small Evelia Rosa Garc\'{\i}a Barroso\\
Departamento de Matem\'atica Fundamental\\
Facultad de Matem\'aticas, Universidad de La Laguna\\
38271 La Laguna, Tenerife, Espa\~na\\
e-mail: ergarcia@ull.es}

\medskip

\noindent {\small Arkadiusz P\l oski\\
Department of Mathematics\\
Technical University \\
Al. 1000 L PP7\\
25-314 Kielce, Poland\\
e-mail: matap@tu.kielce.pl}

\end{document}